\documentclass[3p]{elsarticle}
\usepackage{booktabs}
\usepackage{graphicx}
\usepackage{amsmath}
\usepackage{amssymb}
\usepackage{amsthm}
\usepackage{mathrsfs}
\usepackage{lineno}

\newdefinition{rmk}{Remark}
\newproof{pf}{Proof}
\newproof{pot}{Proof of Theorem \ref{thm2}}
\begin{document}
\begin{frontmatter}
\title{A fast parallel algorithm for solving block-tridiagonal systems of linear equations including the domain decomposition method}

\author{Andrew V. Terekhov}
\ead{andrew.terekhov@mail.ru}
\address{Institute of
Computational Mathematics and Mathematical Geophysics,
630090,Novosibirsk,Russia}
\address{Novosibirsk State University, 630090, Novosibirsk, Russia}
\begin{abstract}
In this study, we develop a new parallel algorithm for solving systems of linear algebraic equations with the same block-tridiagonal matrix but with different right-hand sides. The method is a generalization of the parallel dichotomy algorithm for solving systems of linear equations with tridiagonal matrices \cite{terekhov:Dichotomy}. Using this approach, we propose a parallel realization of the domain decomposition method (\mbox{the Schur} complement method). The calculation of acoustic wave fields using the spectral-difference technique improves the efficiency of the parallel algorithms. A near-linear dependence of the speedup with the number of processors is attained using both several and several thousands of processors. This study is innovative because the parallel algorithm developed for solving block-tridiagonal systems of equations is an effective and simple set of procedures for solving engineering tasks on a supercomputer.
\end{abstract}
\end{frontmatter}
\section{Introduction}
Solving systems of linear algebraic equations (SLAEs) is one of the main problems of computational mathematics. With the advent of multiprocessor computer systems it appeared possible to reduce to some extent computer costs. However in the course of investigation it became evident that most of efficient numerical methods cannot be effectively implemented for supercomputers with many processors.  As supercomputer performance is mainly increased  at the cost of the join of a large number of processors, there arises a necessity to develop the new parallel numerical algorithms for solving SLAEs.

When implementing many numerical techniques, it is required to solve SLAEs with block-tridiagonal matrices\cite{Samarski_Nikolaev,Golub1989,Marsh1981}

\begin{equation}
\label{main-eq} P\mathbf{X}=\left(
\begin{array}{ccccc}
  C_1 & -B_1 &  &   &  \Large 0\\
  -A_2 & C_2 &-B_2  &  &  \\
    & \ddots & \ddots & \ddots &  \\
    &  & -A_{N-1} & C_{N-1} & -B_{N-1} \\
   0&   &  & -A_{N} & C_{N} \\
\end{array}
\right)\left(
\begin{array}{c}
  \bar{\mathbf{X}}_1\\
  \bar{\mathbf{X}}_2  \\
     \dots   \\
    \bar{\mathbf{X}}_{N-1} \\
   \bar{\mathbf{X}}_N \\
\end{array}
\right)=\left(
\begin{array}{c}
  \bar{\mathbf{F}}_1\\
  \bar{\mathbf{F}}_2  \\
     \dots   \\
    \bar{\mathbf{F}}_{N-1} \\
   \bar{\mathbf{F}}_N \\
\end{array}
\right)=\mathbf{F},
\end{equation}
where $A_j,B_j,C_j \in \mathfrak{R}^{\mathrm{M} \times \mathrm{M}}, \; \mathbf{\bar{X}}_j,\mathbf{\bar{F}}_j \in \mathfrak{R}^{\mathrm{M}}$.

By now, various algorithms for solving problem (\ref{main-eq}) on a multi-processor computer system have been developed \cite{block:Ruggieroa,block:Akimova,block:Hirshman,block:Manyu,block:Mehrmann,Chung1995,Bai2007}. But for a multiple solution of SLAEs with the same matrix using the Dichotomy Algorithm this procedure is possible to offer a parallel algorithm with a higher performance as compared to other approaches. The Dichotomy Algorithm is compatible with other algorithms, however it essentially benefits in terms of the time needed for interprocessor interactions. This results from the fact that when implementing the dichotomy process  on a supercomputer it reduces to the calculation of a sum of a series for distributed data thus essentially decreasing the total computing time\cite{terekhov:Dichotomy}.

First the parallel Dichotomy Algorithm was designed for solving SLAEs with the same tridiagonal matrix but different right-hand sides. In \cite{Terekhov2}, the Dichotomy Algorithm was applied to solving SLAEs with Toeplitz tridiagonal matrices. It was shown that for Toeplitz tridiagonal matrices, SLAEs can be effectively solved both with one and several right-hand sides. In \cite{Terekhov2010,fatab2011}, the Dichotomy Algorithm was applied to implement a spectral-difference method of calculation of acoustic and elastic wave fields. This made it possible to effectively use from $2$ up to $8192$ processors per one calculation and to obtain a highly accurate numerical solution of the dynamic problem of elasticity theory. Thus, all the above bears witness to the fact that the Dichotomy Algorithm for solving SLAEs with tridiagonal matrices is a powerful instrument of the numerical modeling. In this paper we propose the new parallel algorithm based on the Dichotomy Algorithm for solving problem (\ref{main-eq}).

When solving many mechanics problems, algorithms based on the domain decomposition method are widespread\cite{Tarek2008,Toselli:2004:DDM,Quarteroni:1999:DDM}. Such an approach has proved its efficiency for calculations on one-processor computers. However with parallel realization of the domain decomposition method, difficulties emerge due to the necessity of implementing efficient algorithms for solving SLAEs. The fact is, efficient methods are, as a rule, difficult to parallelize. We will show that the numerical procedure developed for solving problem (\ref{main-eq})  will allow the effective use of the domain decomposition method (the Schur complement method) for the simulation of acoustic wave fields with thousands of processors.

\section{The Parallel Dichotomy Algorithm for block-tridiagonal matrices }
\subsection{The central idea}
Introduce the following notations:
\begin{itemize}
    \item
    Denote by  $\left\{A\right\}_{l}^{t}$ the matrix obtained from a matrix  $A$ by throwing off all rows and
    columns with the numbers less than  $l$ or greater than $t$.
    \item Denote by  $\left\{\mathbf{V}\right\}_{l}^{t}$ the subvector obtained from a vector
      $\mathbf{V}$ by throwing off the components with the numbers less than  $l$  or greater than  $t$.
    \item  Denote by $\mathbf{e}^\mathrm{L}=\left(1,0,0,...,0\right)^{\mathrm{T}},\;\mathbf{e}^\mathrm{R}=\left(0,...,0,0,1\right)^{\mathrm{T}}$.
\end{itemize}
Omitting unnecessary details, let us formulate a step of the dichotomy process for dividing system (\ref{main-eq}) into two independent subproblems by calculation of the element $\bar{\mathbf{X}}_K$.

\textbf{Algorithm 1}
\begin{enumerate}

  \item Calculate rows of the matrix $P^{-1}_{i\cdot}$ with numbers, where $i=(K-1)M+1,(K-1)M+2,...,KM$.
  \item Calculate the subvector
$$   \bar{\mathbf{X}}_K=\left(P^{-1}_{(K-1)M+1\cdot}\mathbf{F},P^{-1}_{(K-1)M+2\cdot}\mathbf{F},...,P^{-1}_{KM\cdot}\mathbf{F}\right)^{\mathrm{T}}.$$
  \item Transfer from system (\ref{main-eq}) to two independent subsystems by modifying the right-hand side

\begin{subequations}
\begin{equation}\label{suba}
\left\{P\right\}_{1}^{(K-1)M}\left\{\mathbf{X}\right\}_{1}^{(K-1)M}=\left\{\mathbf{F}\right\}_{1}^{(K-1)M}+\mathbf{e}^{\mathrm{R}}\otimes \left(B_{K-1}\mathbf{\bar{X}}_K\right),\; K>1,
\end{equation}

\begin{equation}\label{subb}
\left\{P\right\}_{KM+1}^{NM}\left\{\mathbf{X}\right\}_{KM+1}^{NM}=\left\{\mathbf{F}\right\}_{KM+1}^{NM}+\mathbf{e}^{\mathrm{L}}\otimes \left(A_{K+1}\mathbf{\bar{X}}_K\right),\; K<N.
\end{equation}
\end{subequations}
\end{enumerate}

Further a similar procedure is applied to independent subproblems (\ref{suba}) and (\ref{subb}). Thus, all the components from the solution vector will be calculated in $\lceil \log_2N \rceil$ steps.
Rows of the inverse matrix are stored in the course of calculation and are not recalculated for each right-hand side. As a result, the Dichotomy Algorithm allows "multiplication" of a vector of the right-hand side by the matrix $P^{-1}$ in $O(M^2N\log_2N)$ arithmetical operations, while the direct multiplication would demand $O(M^2N^2)$ operations.  The number of arithmetical operations is decreased because when multiplying a vector by the matrix $P^{-1}$ the information about the structure of the matrix $P$ is used in the Dichotomy Algorithm.

At this point, a consideration of the solution of problem (\ref{main-eq}) with the help of the Dichotomy Algorithm could be completed if there were no an essential complication: it is required to carry out $O(M^3N)$ arithmetical operations as preliminary to the Dichotomy Algorithm in order to calculate rows of the matrix $P^{-1}$\cite{Samarski_Nikolaev}. Such arithmetical costs for large $M$ and $N$ can be unacceptable. Moreover, each processor will require $3M^2N$ RAM cells for storing a copy of the matrix $P$. The use of a supercomputer suggests the solutions of SLAEs of high orders, therefore it is necessary to decrease the required volume of RAM and to minimize the time of preliminary calculations, otherwise it will be impossible to use the Dichotomy Algorithm.
\subsection{An Improved version of the algorithm}
In Algorithm $1$, the basic idea of dividing SLAEs with block-tridiagonal matrices is considered. If the number of processors exceeds the order of the matrix, then auxiliary values $\beta^{\mathrm{R,L}}$, $\mathbf{Z}^{\mathrm{R,L}}$  are introduced\cite{terekhov:Dichotomy,Terekhov2}. But as was noted above, such an approach for block-tridiagonal systems requires high computer costs because its implementation requires solving the original equations system on each processor. Let us explain how to overcome this difficulty.

In \cite{Konovalov,Wang2}, a parallel algorithm based on the superposition principle for solving tridiagonal SLAEs, is proposed. Its central idea is in that the original SLAE reduces to a system of linear equations with a tridiagonal matrix of order $p$, where $p$ is the number of processors. In order to calculate the matrix with the reduced system of equations, on each processor it is necessary to preliminarily solve local subsystems of $\tilde{N}/p$ equations, where $\tilde{N}$ is dimension of a tridiagonal SLAE. After solving the reduced system of equations all the components of the solution vector are independently calculated on each processor.

\begin{figure}[!h]
\begin{center}
\includegraphics[width=0.6\textwidth]{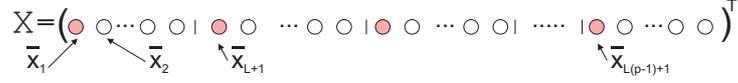}\hfill
\end{center}
\caption{Components of the solution vector to be calculated for dividing the original equations system into subproblems.}
\label{pic:first}
\end{figure}

A similar approach to solving SLAEs with block-tridiagonal matrices is considered in  \cite{block:Akimova}. It consists in the following. The solution to the original system of equations is expressed through $Mp$  of  whilst unknown components from the solution vector (Fig.~\ref{pic:first}):
\begin{equation}
\label{superposition}
\begin{array}{lr}
\mathbf{\bar{X}}_i=\left(\mathbf{U}_i^1\mathbf{U}_i^2...\mathbf{U}_i^\mathrm{M}\right)\mathbf{\bar{X}}_K+\left(\mathbf{V}_i^1\mathbf{V}_i^2...\mathbf{V}_i^\mathrm{M}\right)\mathbf{\bar{X}}_{K+L}+\mathbf{W}_i=U_i\mathbf{\bar{X}}_K+V_i\mathbf{\bar{X}}_{K+L}+\mathbf{W}_i, \\\\   K=1,L+1,2L+1,...,(p-1)L+1;\quad i\in \left[K,K+L\right);\quad L=N/p,\\\\
\mathbf{\bar{X}}_{N+1}=0,
\end{array}
\end{equation} where the matrices $U_i,V_i\in \mathfrak{R}^{\mathrm{M\times M}}$ and the vector $\mathbf{W}_i$  are defined from the solution to subproblems
\begin{subequations}
\begin{equation}\left\{
\begin{array}{llr}
-A_i \mathbf{U}^1_{i-1}+C_i \mathbf{U}^1_{i}-B_{i}\mathbf{U}_{i+1}^1=0,&\mathbf{U}_K^1=\mathbf{e}_1,&\mathbf{U}^1_{K+L}=\mathbf{0},\\
.\ .\ .\ .\ .\ .\ .\ .\ .\ .\ . \ . \ . \ . \ . \ . \ . \ . \ . \ .& \\
-A_i \mathbf{U}^\mathrm{M}_{i-1}+C_i \mathbf{U}^\mathrm{M}_{i}-B_{i}\mathbf{U}_{i+1}^\mathrm{M}=0,&\mathbf{U}_K^\mathrm{M}=\mathbf{e}_M,&\mathbf{U}^\mathrm{M}_{K+L}=\mathbf{0},
\end{array}\right.
\label{sub1}
\end{equation}

\begin{equation}\left\{
\begin{array}{lll}
-A_i \mathbf{V}^1_{i-1}+C_i \mathbf{V}^1_{i}-B_{i}\mathbf{V}_{i+1}^1=0,&\mathbf{V}_K^1=\mathbf{0},& \mathbf{V}^1_{K+L}=\mathbf{e}_1,\\
.\ .\ .\ .\ .\ .\ .\ .\ .\ .\ . \ . \ . \ . \ . \ . \ . \ . \ . \ .& \\
-A_i \mathbf{V}^\mathrm{M}_{i-1}+C_i \mathbf{V}^\mathrm{M}_{i}-B_{i}\mathbf{V}_{i+1}^\mathrm{M}=0,&\mathbf{V}_K^\mathrm{M}=\mathbf{0},& \mathbf{V}^\mathrm{M}_{K+L}=\mathbf{e}_M,
\end{array}\right.
\label{sub2}
\end{equation}
\end{subequations}

\begin{equation}
-A_i \mathbf{W}_{i-1}+C_i \mathbf{W}_{i}-B_{i}\mathbf{W}_{i+1}=\bar{\mathbf{F}}_i,\quad \mathbf{W}_K=\mathbf{0},\quad  \mathbf{W}_{K+L}=\mathbf{0},
\label{sub3}
\end{equation}
where $\mathbf{e}_n$ is a unit vector in the space $\mathfrak{R}^{\mathrm{M}}$.

From (\ref{main-eq}),(\ref{superposition}) we obtain that values of the components $\mathbf{\bar{X}}_K,\; K=1,2L+1,3L+1,...,(p-1)L+1$ can be determined from solving a three-point system of vector equations.

\begin{equation}
\left\{
\begin{array}{lr}
-\left[A_KU_{K-1}\right]\mathbf{\bar{X}}_{K-L}+\left[C_K-A_KV_{K-1}-B_{K}U_{K+1}\right]\mathbf{\bar{X}}_K-\left[B_KV_{K+1}\right]\mathbf{\bar{X}}_{K+L}
=\\\\=\mathbf{F}_K+A_K\mathbf{W}_{K-1}+B_K\mathbf{W}_{K+1},\quad\quad K=1,L+1,2L+1...,(p-1)L+1,\\\\
U_0=V_0=\mathbf{0}.
\end{array}\right.
\label{red1}
\end{equation}

Denote system (\ref{red1}) as $\tilde{P}\tilde{\mathbf{X}}=\tilde{\mathbf{F}}$. For solving system (\ref{red1}) the Dichotomy Algorithm can be applied more effectively than for solving system (\ref{main-eq}). This is due to the fact that reduced system (\ref{red1}) has the dimension $Mp$, while that of original problem (\ref{main-eq}) is $MN$, where $N>p$. As a result, less computer time is needed for a preliminary to the Dichotomy Algorithm as well as a lesser RAM volume $(3M^2N$ vs. $3M^2p)$. Thus, instead of Algorithm $1$ one should use the following algorithm:
\\\\
\textbf{Algorithm 2.}

\begin{enumerate}
  \item The preliminary computations is carried out once  for all the right-hand sides.
  \subitem 1.1 Solve subproblems (\ref{sub1}),(\ref{sub2}) independently on each processor.
  \subitem 1.2 Calculate entries of the matrix $\tilde{P}$ from (\ref{red1}) and send them to all the processors.
  \subitem 1.3 On each processor calculate the required rows of the matrix $\tilde{P}^{-1}$ from $(\ref{red1})$ (for Algorithm $1$).
  \item The stage of calculating solutions is carried out for each right-hand side.
  \subitem 2.1 On each processor solve independently subsystem  (\ref{sub3}).
  \subitem 2.2 Solve system (\ref{red1}) by means of the Dichotomy Algorithm (Algorithm 1).
  \subitem 2.3 In line with (\ref{superposition}) calculate all the components of the solution vector.
\end{enumerate}

At the preliminary step to Algorithm $2$ it is required to solve subsystems $(\ref{sub1}),(\ref{sub2})$. At this stage computer costs are about $O\left(M^3N/p\right)$ arithmetical operations. In order to solve system (\ref{red1}), Algorithm $1$ is used. Therefore it is needed to carry out $O\left(M^3p\right)$ arithmetical operations for calculation of necessary rows of the matrix $\tilde{P}^{-1}$. As entries of the matrix $\tilde{P}$ are distributed among different processors,  the calculation of required rows of the matrix $\tilde{P}^{-1}$ will require interprocessor interactions. The time needed for interprocessor interactions for distributing copies of the matrix $\tilde{P}$ among all the processors will be \footnote{Depending on the algorithm of distribution.} \cite{MPI2}

$$\mathrm{T}^1_{\mathrm{{comm}}}=\alpha \log_2p+\frac{p}{p-1}\beta M^2,$$ $\alpha$--latency, $\beta$--transfer time per byte.

At the second stage of Algorithm $2$, computer costs for solving system (\ref{sub3}) and implementing (\ref{superposition}) will be about $O\left(M^2N/p\right)$. Here the matrix of system (\ref{sub3}) is assumed to be pre-factorized, and the matrices $U_i,V_i$, were computed at the preliminary step. Computer costs of Algorithm $1$ for solving equation (\ref{red1}) are equal to $O\left(M^2\log_2(p)\right)$. Communication costs at the stage of solution calculation are conditioned by the dichotomy process and are estimated as \cite{terekhov:Dichotomy}
$$\quad\mathrm{T}_{\mathrm{comm}}^2\approx\alpha\log^2_2(p)+4M^2\log_{2}(p)\beta.$$

In addition let us note that the necessary volume of RAM at this stage will be $O\left(M^2N/p+M^2\log_2p\right)$, while at the preliminary step it makes $O\left(M^2N/p+M^2p\right)$.

\subsection{Numerical experiments}
Let us consider the problem of solving a system of linear equations of the form of (\ref{main-eq}) with dimensions of blocks $M=60,150$ for $N$ from $2048$ up to $65536$.

\begin{figure}[!h]
\begin{center}
\includegraphics[width=0.5\textwidth]{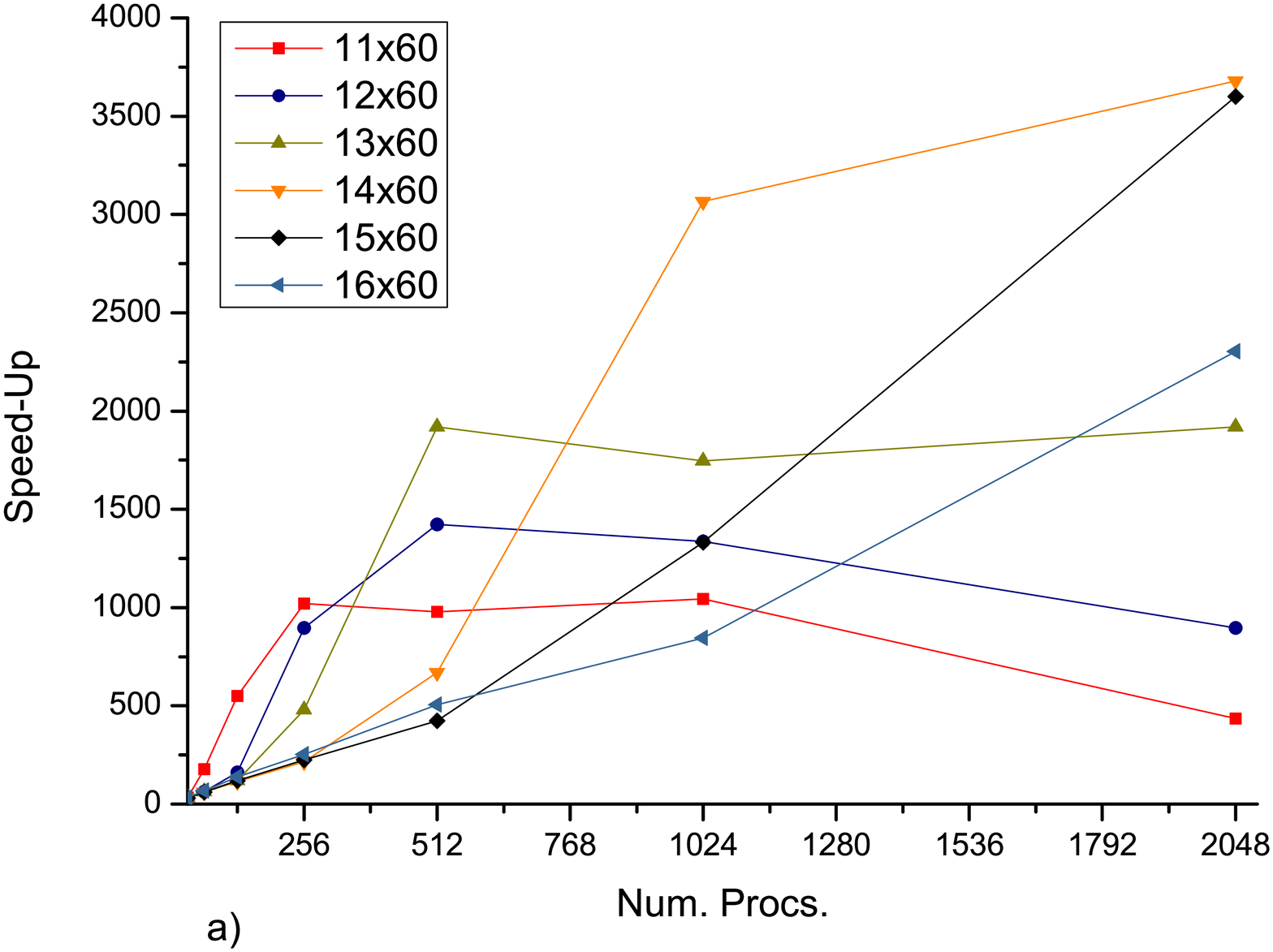}\hfill
\includegraphics[width=0.5\textwidth]{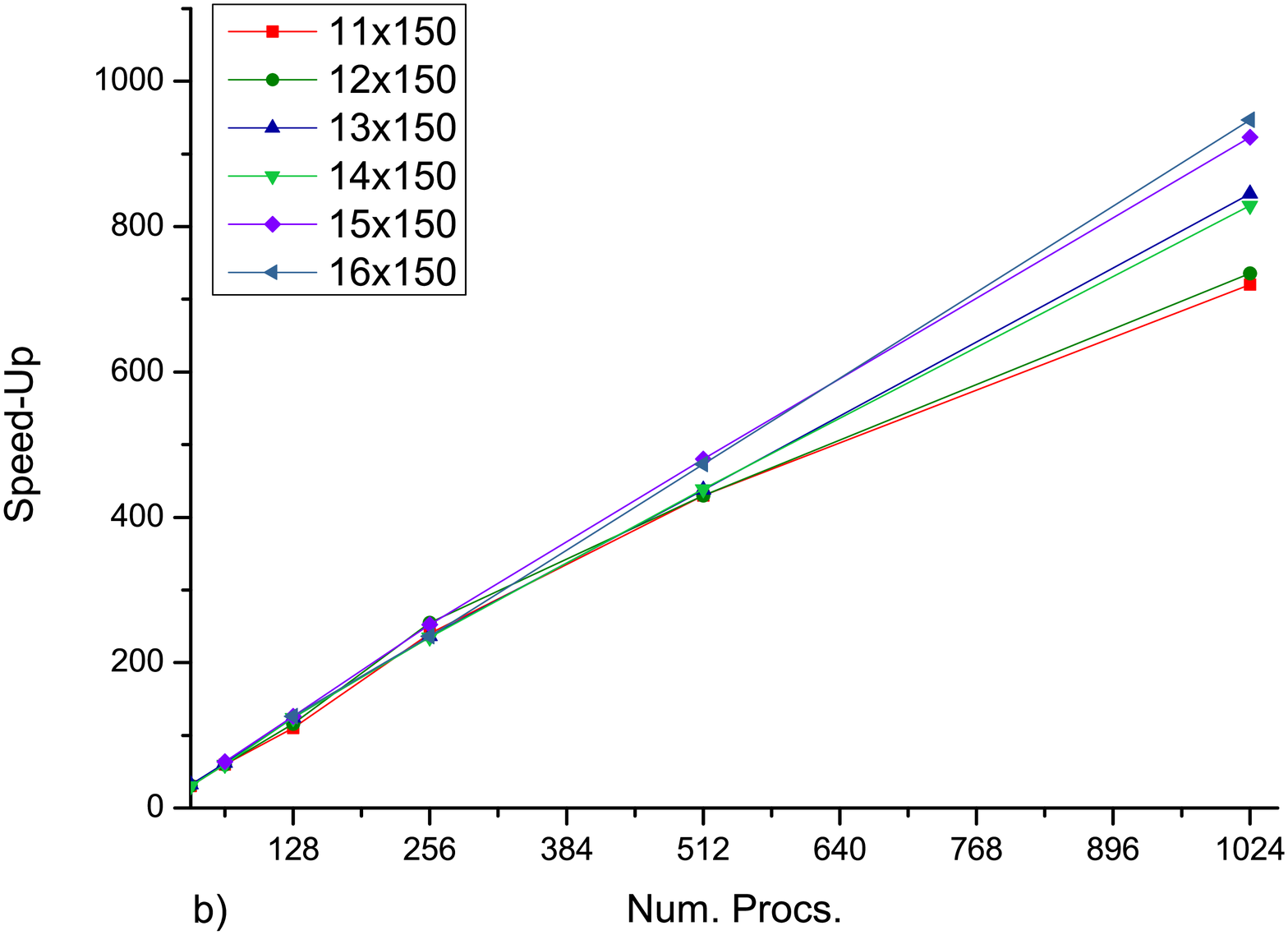}
\end{center}
\caption{Dependence of the speedup value on the number of processors for $M=60$ (a) and $M=150$ (b) on various $N$.}
\label{pic:m150}
\end{figure}

Numerical procedures were implemented in Fortran-90 using MPI library, calculation being performed on "MBC-100k"\ supercomputer of the Interdepartment Supercomputer Center of the Russian Academy
of Sciences, (the 62-nd position in Top-500\cite{top500}, November 2010). The results of the experiments conducted are given in Tables~\ref{table1},\ref{table2} and in Fig.\ref{pic:m150}.

\begin{table}[!h]
\center \small
\begin{tabular}{lllllllllllll}

  \hline
  $N \times M$ & \multicolumn{2}{c}{$2^{11}\times 60$}& \multicolumn{2}{c}{$ 2^{11}\times 150$}& \multicolumn{2}{c}{$2^{12}\times 60$}& \multicolumn{2}{c}{$ 2^{12}\times 150$}&\multicolumn{2}{c}{$2^{13} \times 60$}&\multicolumn{2}{c}{$2^{13} \times 150$} \\ \cmidrule(r){2-3} \cmidrule(r){4-5} \cmidrule(r){6-7} \cmidrule(r){8-9}  \cmidrule(r){10-11} \cmidrule(r){12-13}
   NP &  $ {\mathrm{Pre}}$&  $ \mathrm{Exe}$  & $ \mathrm{Pre}$&$ \mathrm{Exe}$&$ \mathrm{Pre}$&$ \mathrm{Exe}$& $\mathrm{Pre}$&$ \mathrm{Exe}$& $\mathrm{Pre}$&$ \mathrm{Exe}$ & $\mathrm{Pre}$&$ \mathrm{Exe}$\\ \hline
  16 &1.8&3e-2&30&0.18 &3.5&5.6e-2&56&0.35&7&0.12&115&0.74 \\
  32   &1&1.5e-2&17&9.4e-2&1.9&2.9e-2&32&0.18&3.8&6e-2&60&0.37 \\
  64   &1.32&2.7e-3&15&4.8e-2&2.7&1.4e-2&22&9.4e-2&3&3e-2&36&0.19 \\
  128  &2&8.7e-4&19&2.6e-2&2.1&5.5e-3&22&4.8e-2&3&1.6e-2&30&9.5e-2   \\
  256  &3.2&4.7e-4&32&1.2e-2&3.2&1e-3&34.7&2.6e-3&3.36&4e-3&38&5e-2 \\
  512  &5.6&4.9e-4&62&6.7e-3&5.7&6.3e-4&61&1.3e-2&5.6&1e-3&65&2.7e-2  \\
  1024 &10&4.6e-4&124&4e-3&9.9&6.7e-4&126&7.6e-3&9.9&1.1e-3&123&1.4e-2  \\
  2048 &48.8&1.1e-3&$\Theta$&$\Theta$&19.44&1.3e-4&$\Theta$&$\Theta$&47&1e-3& $\Theta$&$\Theta$ \\\hline
\end{tabular}
\caption{Preliminary time ($\mathrm{\mathbf{Pre}}$) and execution time (${\mathrm{\mathbf{Exe}}}$).$2^{11}=2048$.}
\label{table1}
\end{table}

\begin{table}[!h]
\center \small
\begin{tabular}{lllllllllllll}

  \hline
    $N \times M$ & \multicolumn{2}{c}{$2^{14}\times 60$}& \multicolumn{2}{c}{$ 2^{14}\times 150$}& \multicolumn{2}{c}{$2^{15}\times 60$}& \multicolumn{2}{c}{$ 2^{15}\times 150$}&\multicolumn{2}{c}{$2^{16} \times 60$}&\multicolumn{2}{c}{$2^{16} \times 150$} \\ \cmidrule(r){2-3} \cmidrule(r){4-5} \cmidrule(r){6-7} \cmidrule(r){8-9}  \cmidrule(r){10-11} \cmidrule(r){12-13}
   NP &  $ {\mathrm{Pre}}$&  $ \mathrm{Exe}$  & $ \mathrm{Pre}$&$ \mathrm{Exe}$&$ \mathrm{Pre}$&$ \mathrm{Exe}$& $\mathrm{Pre}$&$ \mathrm{Exe}$& $\mathrm{Pre}$&$ \mathrm{Exe}$ & $\mathrm{Pre}$&$ \mathrm{Exe}$\\ \hline
  16  &15&0.23&220&1.4&28&0.45&$\Theta$&$\Theta$&60&0.95&$\Theta$&$\Theta$\\
  32  &7.5&0.12& 117&0.73&15&0.24&231&1.5&29&0.47&$\Theta$&$\Theta$\\
  64  &4.3&6e-2& 64&0.37&8.75&0.12&121&0.74&15 &0.22&235&1.48\\
  128 &4.2&3.2e-2&44&0.18&5.15&6e-2& 72&0.38&8.5&0.11&129&0.75 \\
  256 &4.1&1.7e-2&45&9.5e-2&5&3.2e-2&60 &0.19&6.9&6e-2&89&0.4\\
  512 &5.6&5.5e-3&68&5.1e-2&14.6&3.6e-2&75&0.1& 7.0 &3e-2&89&0.2 \\
  1024&14.4&8.8e-4&126&2.7e-2&25&5.4e-2&130&5.5e-2&10.3&1.8e-2&136&0.1 \\
  2048&47&1e-3&$\Theta$&$\Theta$&48&2e-3&$\Theta$&$\Theta$&48&6.6e-3 &$\Theta$&$\Theta$\\\hline
\end{tabular}
\caption{Preliminary time ($\mathrm{\mathbf{Pre}}$) and execution time (${\mathrm{\mathbf{Exe}}}$).$2^{14}=16384$.}
\label{table2}
\end{table}

Based on the data obtained let us note the following:
\begin{itemize}
\item In all the test calculations, the value of dependence of the speedup value on the number of processors was near linear.

\item For matrices with $M=60$-blocks, starting with a certain $p>p_0$ the dependence of the speedup value on the number of processors was superlinear. This is due to increasing a general number of processors which, in turn, causes a decrease in the data volume of the problem for one processor. Thus, this allows a more effective use of a high-speed cache memory. A similar effect was achieved with a parallel realization of the ADI method\cite{terekhov:Dichotomy}.

\item The preliminary time depends on the number of processors used. With a minor amount of processors, the main costs fall on solving problems (\ref{sub1}),(\ref{sub2}) and decrease with the growth of the number of processors. But starting with a certain $p>p_0$, preliminary costs for Algorithm $1$ that are required for a subsequent solution of problem (\ref{red1}) become dominating.

\item For matrices with $M=150$-blocks and the number of processors $p=2048$ it appeared impossible to carry out preliminaries in a reasonable time. This is because the matrix of the reduced system has $M^2p$ dimension and with $p=2048$ processors cannot be completely located in RAM of one computer unit. The use of disk memory has considerably decreased the performance.

\item For the parameters $M=150,\ N=2^{16}$, $p=16,32$ and $M=150,\ N=2^{15}$, $p=16$ the insufficient volume of RAM because of a small number of processors did not allow solving problem (\ref{sub1}),(\ref{sub2}) in a reasonable time.

\end{itemize}

The numerical experiments have shown that the Dichotomy Algorithm provides a high efficiency of using supercomputer resources. When realizing the Dichotomy Algorithm in terms of numerical procedures one should pay attention to available volume of RAM, because as compared to iterative techniques both for the Dichotomy Algorithm and for most of direct methods of solving SLAEs a larger volume of RAM is needed.

\section{Acoustic Solver}
To gain greater insight into the Dichotomy Algorithm efficiency for solving applied problems of numerical modeling, in the cylindrical coordinate system $(r,z)$, in the half-space $z\geq 0$ we will consider the problem of modeling the propagation of acoustic waves from a point source

\begin{equation}
\begin{array}{llr}
\displaystyle {\rho({\bf x} )}\frac{\partial^2 p}{\partial
t^2}({\bf x},t)=\nabla \left[\kappa({\bf x})\, \nabla p({\bf
x},t)\right]+\frac{1}{2\pi}\frac{\delta({\bf x-x_0})}{r}f(t),&
t>0,\quad {\bf x}=(r,z),
\end{array}
\label{acoustic_problem}
\end{equation}
{ where $p(\mathbf{x},t)$ is the acoustic pressure, $\rho(\mathbf{x})$ is the density perturbations, $\sqrt{{\kappa(\mathbf{x})}/{\rho(\mathbf{x})}}$ is the sound velocity, $\mathbf{x}_0$ is the source coordinates.}
Suppose that problem (\ref{acoustic_problem}) is solved with homogeneous initial conditions.

A parallel version of the spectral-difference method for solving (\ref{acoustic_problem}) was considered in \cite{Terekhov2010,fatab2011}. The Laplace operator was selected as preconditioning operator. This allowed us to provide a high rate of convergence for media with moderate contrast. The use of the Dichotomy Algorithm for solving tridiagonal SLAEs made possible to attain a high calculation rate. However when a medium model includes zones of high and relatively low velocities, using the Laplace operator as preconditioning does not provide a high convergence rate of the iterative process for solving SLAEs. If it appears possible to distinguish macro-zones in the medium model, where the sound speed is constant or is slightly diverse, then it makes sense to use the domain decomposition Method. Parallel versions of the domain decomposition method were proposed  rather a long time ago \cite{Smith:1996,Bitzarakis1997a} and recently algorithms with graphics accelerators have been offered\cite{DD:GPU}. In this paper, the domain decomposition method based on the Schur complement will be used for decreasing the number of arithmetical operations for solving difference equations but not as parallelization instrument. Thus, in the case under consideration the number of processors used and the number of subdomains will be independent values. The efficiency of using supercomputer resources will be completely provided at the cost of employing the Dichotomy Algorithm.

\begin{figure}[!h]
\begin{center}
\includegraphics[width=0.5\textwidth]{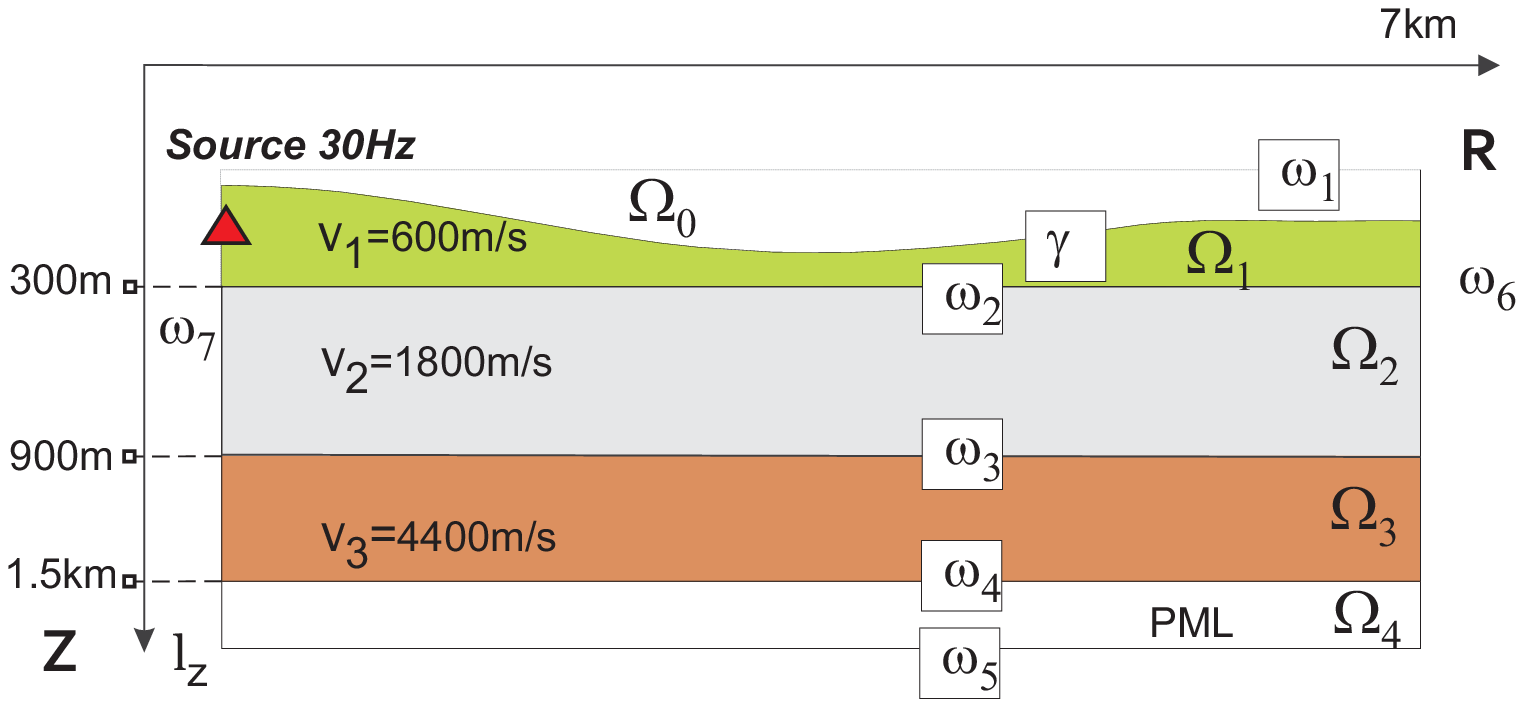}\hfill
\includegraphics[width=0.5\textwidth]{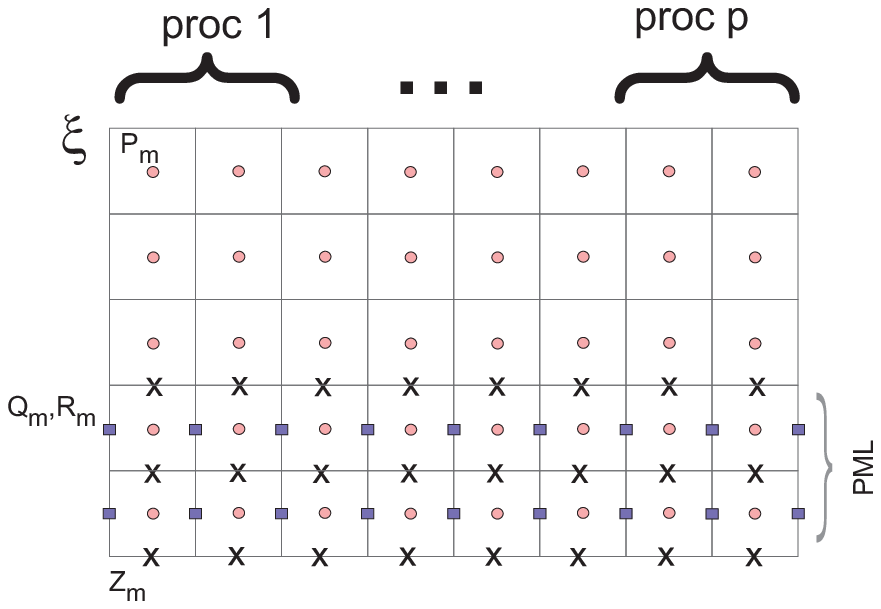}
\end{center}
\caption{Medium model (a) and solution mesh (b).}
\label{pic1:model}
\end{figure}

Figure~\ref{pic1:model}.a. presents a medium model, for which it is reasonable to use the domain decomposition method.
When solving applied geophysics problems it is often required to calculate a wave field with an arbitrary geometry of the free surface \cite{Sheriff1995}, therefore such the medium model will include a relief.

In addition, to exclude non-physical reflections from the fictitious boundary $\omega_4$ in the subdomain $\Omega_4$, the PML absorbing boundary conditions will be realized\cite{Berenger1994,Chew1994}

\begin{equation}
\left\{
\begin{array}{ll}
 \left(\frac{\partial}{\partial t}+\sigma_z\right)p-\rho_0 c^2\left(\frac{1}{r}\frac{\partial}{\partial r }\left(v_r r\right)+\sigma_z\frac{1}{r}\frac{\partial}{\partial r}\left(qr\right)+\frac{\partial v_z}{\partial z}\right)=0,\\\\
 \frac{\partial v_r}{\partial t}-\frac{1}{\rho_0}\frac{\partial }{\partial r} p=0,\quad\quad   \left(\frac{\partial}{\partial t}+\sigma_z\right) v_z-\frac{1}{\rho_0}\frac{\partial }{\partial z} p=0,\quad\quad
\frac{\partial q}{\partial t}=v_r,
\end{array}
\right.
\label{pml_eq}
\end{equation}
where the absorbing layers profile is given by the function $\sigma_z(z)=\frac{(\nu+1)c_p}{2L_{\mathrm{PML}}}\log\left(\frac{1}{|\chi|}\right)\left[\frac{(z-z_0)}{L_{\mathrm{PML}}}\right]^\nu$, $\chi$ is a user-tunable reflection coefficient, $\nu$ is the degree of the polynomial attenuation, $c_p$ is the wave velocity, $L_{\mathrm{PML}}$ is a width of PML region.
\subsection{The Laguerre transform}
Let us seek for a solution to problem (\ref{acoustic_problem}) as a Fourier series in the Laguerre functions \cite{Mikhailenko1999}

\begin{equation}
p(\mathbf{x},t)=(\eta t)^{\frac{\alpha}{2}}\sum_{m=0}^{\infty}P_m({\bf
x})l^{\alpha}_m(\eta t), \label{series_lag1}
\end{equation}
where $l^{\alpha}_m(\eta t)$ are the orthonormal Laguerre functions\cite{abramowitz+stegun}, $m$ is Laguerre polynomial degree, $\alpha$ is the order of Laguerre functions and $\eta$ is the transformation parameter.
Applying the Laguerre transform to (\ref{acoustic_problem}), we obtain a series of problems for defining the expansion factors

\begin{equation}
\left\{
\begin{array}{ll}
\label{laguerre_h}   \nabla \left[\kappa({\bf x})\,
\nabla P_m({\bf x})\right]-\rho({\bf x})\frac{\eta^2}{4}P_m({\bf
x})= -\frac{1}{2\pi}\frac{\delta({\bf
x-x_0})}{r}f_m+\rho({\bf x})\eta^2\sqrt{\frac{m!}{(m+\alpha)!}}\sum_{k=0}^{m-1}(m-k)\sqrt{\frac{(k+\alpha)!}{k!}}P_k({\bf x})\ \mathrm{in} \ \cup_{i=1}^3 \Omega_i,\\
  \frac{\partial P_m}{\partial n}=0 \  \mathrm{on} \  \gamma, \quad \frac{\partial P_m}{\partial r}=0 \  \mathrm{on} \  \omega_7, \quad
  P_m=0\  \mathrm{on} \  \omega_6,
\end{array}
\right.
\end{equation}
где $f_m=\int_0^{\infty}f(t)(\eta t)^{-\frac{\alpha}{2}}l^{\alpha}_{m}(\eta t)dt$.

Applying the Laguerre transform to equations (\ref{pml_eq}) and introducing the notation
$$\Phi(X_m)\equiv \eta\sqrt{\frac{(m+1)!}{(m+\alpha+1)!}}\sum_{k=0}^{m}\sqrt{\frac{(k+\alpha)!}{k!}}X_k,$$ we obtain the following system of equations

\begin{equation}
\left\{
\begin{array}{ll}
\left(\left[1+\frac{2\sigma_z}{\eta}\right]\frac{1}{r}\frac{\partial}{\partial r }\left(R_m r\right)+\frac{\partial Z_m}{\partial z}\right)-\frac{1}{\rho_0c^2}\left(\frac{\eta}{2}+\sigma_z\right)P_m=\Phi(P_{m-1})+\frac{2\sigma_z}{r\eta}\frac{\partial}{\partial r}(r\Phi(Q_{m-1})) ,\\\\
 \frac{1}{\rho_0}\frac{\partial P_m }{\partial r } -\frac{\eta}{2}R_m=\Phi(R_{m-1}),\quad\quad
 \frac{1}{\rho_0}\frac{\partial P_m }{\partial z } -\left(\frac{\eta}{2}+\sigma_z \right)Z_m=\Phi(Z_{m-1}),\quad
 \frac{\eta}{2}Q_m+\Phi(Q_{m-1})=R_m.
\end{array}\right.
\label{laguerre_h2}
\end{equation}
Here $R_m,Z_m,P_m,Q_m$ are expansion factors in the Fourier-Laguerre series for the functions $v_r,v_z,p,q$.

\subsection{Domain decomposition}
In the domain $\Omega=\bigcup_{i=0}^4\Omega_i$ (Fig.~\ref{pic1:model}.a) and $\Omega_i\bigcap\Omega_j=\emptyset$ when $i\neq j$ introduce a rectangular mesh $\xi$ (Fig.~\ref{pic1:model}.b). Inside the domains $\Omega_i,\; i=1,2,3$ on the mesh $\xi$ approximate problem (\ref{laguerre_h}), and in the subdomain $\Omega_4$ approximate equation (\ref{laguerre_h2}) for the PML absorbing boundary conditions.

As the approximation of elliptic equations is widely covered in the literature \cite{Samarskii2001,Strikwerda2004,FEM2}, we will only mention that for solving equation (\ref{laguerre_h}) a five-point scheme of second order of accuracy was used that was constructed by the finite volume method. For solving equation (\ref{laguerre_h2}), we made use of the scheme of second order of accuracy on the staggered mesh(Fig.~\ref{pic1:model}.b).

To reduce the dependence of the number of arithmetical operations on the contrast of the medium, let us dwell on the domain decomposition method based on the Schur decomposition \cite{Tarek2008,Toselli:2004:DDM}. To this end the mesh nodes are enumerated in the following order: first the nodes from $\Omega_1,\Omega_2,\Omega_3,\Omega_4$, and then those belonging to the boundaries $\omega_2,\omega_3,\omega_4$. Then the difference problem for equations (\ref{laguerre_h}),(\ref{laguerre_h2}) is written down as SLAE\cite{Tarek2008,Toselli:2004:DDM}

\begin{equation}
\left[
\begin{array}{ccccc}
A_{11} &0&0&0&A_{1\Gamma}\\
0 &A_{22}&0&0&A_{2\Gamma}\\
0 &0&A_{33}&0&A_{3\Gamma}\\
0 &0&0&A_{44}&A_{4\Gamma}\\
A^{\mathrm{T}}_{1\Gamma}&A^{\mathrm{T}}_{2\Gamma}&A^{\mathrm{T}}_{3\Gamma}&A^{\mathrm{T}}_{4\Gamma}&A_{\Gamma\Gamma}
\end{array}\right]\left\{
\begin{array}{ccccc}
\mathbf{x}_1\\
\mathbf{x}_2\\
\mathbf{x}_3\\
\mathbf{x}_4\\
\mathbf{x}_{\Gamma}
\end{array}\right\}=\left\{
\begin{array}{ccccc}
\mathbf{f}_1\\
\mathbf{f}_2\\
\mathbf{f}_3\\
\mathbf{f}_4\\
\mathbf{f}_{\Gamma}
\end{array}\right\}
,
\label{linear_decomp}
\end{equation}
where each $\mathbf{x}_i$ represents the subvector of unknowns that are interior to subdomain $\Omega_i$ and $x_{\Gamma}$ represents the vector of all interface unknowns.

The matrix $A_{jj}$ corresponds to the difference problem for equation (\ref{laguerre_h}) in the interior of the subdomain $\Omega_j,\;j=1,2,3$, while the matrix $A_{44}$ -- for equation (\ref{laguerre_h2}).

First let us calculate components belonging to the boundaries $\omega_{1,2,3}$. To this end we solve the system of equations
\begin{equation}
\begin{array}{ll}
S\mathbf{x}_{\Gamma}=\mathbf{f}_{\Gamma}-\sum_{i=1}^{4}A_{i\Gamma}^{\mathrm{T}}A_{ii}^{\mathrm{-1}}\mathbf{f}_{i},
\end{array}
\label{schur_decomp}
\end{equation}
where the Schur complement $S$ is defined by $S=A_{\Gamma\Gamma}-\sum_{i=1}^{4}A^{\mathrm{T}}_{i\Gamma}A^{\mathrm{-1}}_{ii}A_{i\Gamma}$.

Once $x_{\Gamma}$ is determined, the complete solution in the interior of the subdomains is obtained from
\begin{equation}
\mathbf{x}_{i}=A_{ii}^{\mathrm{-1}}\left(\mathbf{f}_{i}-A_{i\Gamma}\mathbf{x}_{\Gamma}\right),\quad \text{for}\; i=1,2,3,4.
\label{domains_equation}
\end{equation}

For the matrix $S$ be calculated not in the explicit form, we use the conjugate gradient method(the CG method)\cite{Saad} for solving problem (\ref{schur_decomp}). To implement the CG method, it is necessary to solve the two problems. The first one is in that multiplication of a vector by the matrix $S$ requires parallelization of efficient procedures for the multiple inversion of the matrices $A_{ii}$. In addition, the matrix $S$ is ill-conditioned, hence it is required to use the preconditioning procedure. Further we will show that for solving such subproblems, one can build efficient parallel procedures based on the Dichotomy Algorithm.

\subsubsection{Multiplication of a vector by the matrix $S$}
It is evident that the main computer costs are required for the multiple inversion of the matrices $A_{ii}\;, i=1,2,3,4$. These problems are considered to be uniformly distributed among $p$ processors according to Fig.~\ref{pic1:model}.b. Consider parallel procedures for multiplication of a vector by the matrices $\left(A^{\mathrm{T}}_{i\Gamma}A^{\mathrm{-1}}_{ii}A_{i\Gamma}\right),\; i=1,2,3,4$.

\textbf{a. Solution to elliptic equations in the subdomain $\Omega_1$.}
The difference problem for equation (\ref{laguerre_h}) in the subdomain $\Omega_1$ is in agreement with a system of linear algebraic equations with the matrix $A_{11}$. An arbitrary geometry of the free surface $\gamma$ can be taken into account in different ways: irregular grids, the method of Lagrange multipliers, the method of fictitious domains \cite{Glowinski2007,Saulev1963,Ramiere}, conformal mapping \cite{Comformal_Mapping}.
In \cite{fatab2011} it was shown that for calculation of wave fields for long durations of time one should use grids with a high spatial resolution $h_{r,z}\approx1/200 \lambda_{min}\div 1/100 \lambda_{min}$, where $\lambda_{min}$ is a minimum wavelength. If the sound velocity close to the free surface is not high, then due to necessity of using a small mesh size the free boundary $\gamma$ can be smoothed along the boundaries of the nearest cells. In practice, an admissible error for defining the depth of layers bedding, for example for the West Siberia region, makes up about several meters, that is why the mesh size equal to a few centimeters allows approximating with a sufficient accuracy the relief. Such an approximation makes possible to carry out calculations sufficiently fast, which is more reasonable in terms of efficiency. To make use of the approach in question, we apply an algebraic version of the method of fictitious domains, that is the fictitious components technique \cite{Astrakhantsev1978,Marchuk1982}, whose idea is in that a subvector $\phi_0$ being the solution to SLAE with a positive semi-definite matrix
\begin{equation}
\left[
\begin{array}{cc}
A_{11} &0\\
0&0
\end{array}\right]\left\{
\begin{array}{c}
\phi_0\\
\phi_1
\end{array}\right\}=\left\{
\begin{array}{c}
f_0\\
0
\end{array}\right\}
\label{fictitious}
\end{equation}
will also be the solution to the system $A_{11}\phi_{0}=f_0$. Let a matrix $C$ correspond to the difference problem for the operator $L_h\equiv\Delta_h-d^2,\ d\in \mathfrak{R}$ in the subdomain $\Omega_0\cup \Omega_1$. Then system (\ref{fictitious}) can be solved by the GMRES($k$) method with the preconditioning matrix $C$ per the number of iterations independent of the mesh size \cite{Astrakhantsev1978}. In this case, the main macro-operation is in the inversion of the operator $L_h$ thus allowing the use of the Dichotomy Algorithm for the effective parallelization.

\textbf{b. Solution to elliptic equations in the subdomains $\Omega_{2,3}$.} Let matrices $A_{22}$ and $A_{33}$ correspond to the difference problem for equation (\ref{laguerre_h}) for the subdomains $\Omega_2$ and $\Omega_3$, respectively. To multiply a vector by the matrices $A_{22}^{-1},A_{33}^{-1}$ , it is possible to use the method of separation of variables\cite{FFT} with arithmetical operations costs $O(N\log N)$, where $N$ is a common number of mesh nodes in the subdomain. However it appears possible to calculate the product of a vector by the matrices $A_{i\Gamma}^{\mathrm{T}}A_{ii}^{-1}A_{i\Gamma},\; i=2,3$ with essentially lesser arithmetical costs. For the difference problems $A_{ii}\hat{\mathbf{y}}=A_{i\Gamma}\hat{\mathbf{f}},\; i=2,3$ the mesh function $A_{i\Gamma}\hat{\mathbf{f}}$ takes nonzero values only in boundary nodes of the subdomain. For this right-hand side the calculation of the direct Fourier transform will demand only $O(N)$ arithmetical operations\cite{FFT}. To multiply a vector by the matrix $A_{i\Gamma}^{\mathrm{T}}$, it is sufficient to calculate components from the vector $\hat{\mathbf{y}}$ that correspond to boundary mesh nodes for the subdomain. The inverse Fourier transform for defining the solution only in the boundary mesh nodes can be carried out in $O(N)$ arithmetical operations. Taking into account the fact that the solution to tridiagonal SLAE in terms of the method of separation of variables will demand $O(N)$ arithmetical operations, the final assessment of the number of arithmetical operations for multiplying a vector by the matrix $A_{i\Gamma}^{\mathrm{T}}A_{ii}^{-1}A_{i\Gamma}$ will be $O(N)$.
\label{xxx1}

\textbf{c. Solution to elliptic equations in the subdomain $\Omega_{4}$.}
Let a matrix $A_{44}$ correspond to the difference problem for the PML equations (\ref{laguerre_h2}) in the subdomain $\Omega_4$. The width of a PML region is, as a rule, found within the limits of $20$ up to $50$ mesh nodes, while the number of cells in the radial direction is considerably larger, that is $N_r\gg N_z$. With the above enumeration of unknowns, the matrix $A_{44}$ will be a band matrix of order $3N_zN_r$ with the bandwidth $3N_z$, where factor $3$ is conditioned by the necessity of computing the three components $v_r,v_z,p$ for the PML region. As $N_z\ll N_r$ and with allowance for the ill conditioning of the matrix $A_{44}$, for multiplying a vector by the matrix $A_{44}^{-1}$ it seems reasonable to use the Dichotomy Algorithm (Algorithm $2$) for block-tridiagonal matrices\footnote{For the band matrix, the submatrices $A_i,B_i$ from (\ref{main-eq}) are upper triangular and lower triangular matrices, respectively.}.

Thus, all the procedures of solving the local subproblems include the Dichotomy Algorithm. With allowance for the results of computer experiments from the previous section, one should expect that the dependence of the speedup value  on the number of processors for the multiplication of a vector by the matrix $S$ will be close to the linear one.

\subsubsection{Preconditioning}

By now there have been developed relatively many sequential versions of preconditioning procedures for solving problem (\ref{schur_decomp}) \cite{Tarek2008,Toselli:2004:DDM,chan}. However for supercomputers a class of effective preconditioners is essentially less. In this paper, we use a preconditioner based on the probing technique \cite{Chan1992,Tarek2008,Toselli:2004:DDM}. The operator $S$ is approximated by an operator $B$ on a certain subspace, the latter being constructed so as $B$ be readily invertible. In this case the matrix $B$ will be a band one. The probing technique does not demand the knowledge about the structure of the operator $S$ and is a purely algebraic approach. To calculate the matrix $B$, one should realize the multiple multiplication of the matrix $S$ by specially selected vectors $\mathbf{p}_l,1\leq l\leq 2d+1$, where $d$ is the bandwidth of the matrix $B$. This procedure was discussed in the previous section. As the bandwidth is essentially less than the order of the matrix $B$, it appears possible to use the Dichotomy Algorithm for block-tridiagonal matrices when solving SLAEs with a band matrix $B$.

\subsection{Numerical experiments}
Let the size of the computational domain be $L_r=7km,\;L_z=1.5km$. A point source is located on the symmetry axis at a depth of $15m$ from the free surface; the time dependence being given as

\begin{equation}
\label{functiont}
 f(t)=\exp\left[-\frac{(2\pi
f_0(t-t_0))^2}{g^2}\right]\sin(2\pi f_0(t-t_0)),
\end{equation}
where $f_0=30\mathrm{Hz},\;t_0=0.2s,\;g=4$. The number of addends in series (\ref{series_lag1}) was $n=6000$; the expansion parameters were $\alpha=5,\;\eta=1800$. For the PML boundary conditions, the following parameters were selected:
$L_{\mathrm{PML}}=30h_z$, $c_p=4400m/s$, $\nu=2$, $\chi=10^{-6}$.

The issues concerning the spectral algorithm based on the Laguerre transform were studied in \cite{Terekhov2010,fatab2011}, therefore we will dwell on performance and efficiency of the parallel algorithm.

The matrices $S$ and $B$ are not spectrally equivalent, therefore with decreasing the mesh step the number of iterations of the conjugate gradient method for solving (\ref{schur_decomp}) will increase \cite{Chan1992}. However increasing the bandwidth of the matrix $B$, denoted as $d$ makes possible to decrease the number of iterations (Table~\ref{table33}). With a twofold decrease of the mesh step the value of the parameter $d$ should twofold be increased for the number of iterations of the CG method be not increased. Thus, preconditioning based on the probing technique allows a considerable decrease in computer costs, while the Dichotomy Algorithm makes possible to efficiently solve SLAEs with the preconditioning matrix.

\begin{table}[!h]
\center \small
\begin{tabular}{llllllllll}

  \hline
    $N_z \times N_r$ & \multicolumn{3}{c}{$3046\times 16384$}& \multicolumn{3}{c}{$ 6078 \times 32768$}& \multicolumn{3}{c}{$12134\times 65536$}\\ \cmidrule(r){2-4} \cmidrule(r){5-7} \cmidrule(r){8-10}
    d&$ {\mathrm{Gen}}$ & $\mathrm{Total}$&$ \mathrm{Iter}$&$ {\mathrm{Gen}}$ & $\mathrm{Total}$&$ \mathrm{Iter}$&$ {\mathrm{Gen}}$  & $\mathrm{Total}$&$ \mathrm{Iter}$

    \\ \hline
  no prec  &-&370&4000&-&1550&5600&-&13900&9710 \\
  3  &0.75&2.28&10 &1&13 &6.6&4.51&40&16 \\
  5  &0.57&2&9&1.6&4.11&5.6&8.3&38&14 \\
  7  &1.3&1.9&8&2.14&5.61&10&17&36&13 \\
  11  &1.6&1.2&6&3.43&5&9&17.5&34&11 \\
  21  &2.8&1.13&5&6.61&4.45&7&33&32&9 \\
  51  &6.9&0.88&3 &16.24&3.6&4&87&27&6\\
  101  &-&-&-&34&3.54&3&166&26&4 \\
  \hline
  \end{tabular}
\caption{$(\mathrm{\mathbf{Gen}})$ is the time of calculation of the preconditioning matrix $\mathrm{B}$. $(\mathrm{\mathbf{Total}})$ is the time needed for solution to one problem of the form of (\ref{schur_decomp}). In this case it is necessary to carry out $(\mathrm{\mathbf{Iter}})$ iterations. $(\mathrm{\mathbf{d}})$ is the bandwidth of the preconditioning matrix, the number of processes being constant $p=256$.}
\label{table33}
\end{table}

\begin{figure}[!h]
\begin{center}
\includegraphics[width=1\textwidth]{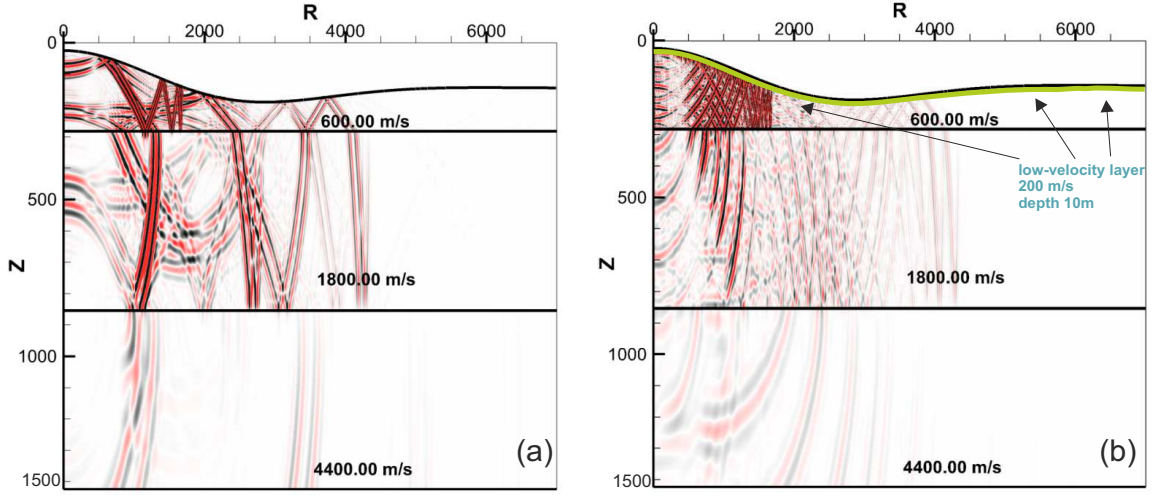}
\end{center}
\caption{Snapshots for the wave field at $t=3s$ (a) for model pic.\ref{pic1:model}.a and same model with additional low-velocity layer (b). $\mathrm{N_z}\times \mathrm{N_r}=12134\times65536$.
}
\label{pic:snapshot}
\end{figure}

A feature of the parallel algorithm proposed is in that before solving a series of problems (\ref{laguerre_h}),(\ref{laguerre_h2}), it is required to conduct preliminary calculations. The time assessments for preliminary calculations $(\mathrm{P})$  are given in Table~\ref{table44}, which also represents the preliminaries for the Dichotomy Algorithm for inverting the matrices $A_{ii},\; i=1,2,3,4$ as well as the costs for the calculation of the preconditioning matrix $B$. The number of terms in series (\ref{series_lag1}) for long time durations makes up several thousands, that is why the time needed for the preliminaries can be neglected. This is because of their smallness as compared to the general computation time $(\mathrm{T})$.

The smaller speedup coefficient $(\mathrm{S})$ (Table~\ref{table44}) as compared to the Poisson equation solution\cite{terekhov:Dichotomy} is due to the necessity of complementary interprocessor communications for the GMRES($k$) method for the matrix $A_{11}$ inversion. Moreover, the multiple inversion of the preconditioning operator $C$ in the interior of the small subdomain $\Omega_0\cup \Omega_1$ causes an increase in the communication time as related to the computation time and, hence, the scalability of the parallel algorithm decreases.

In Section \ref{xxx1} it was shown that the multiplication of a vector by the matrices $A_{i\Gamma}^{\mathrm{T}}A_{ii}^{-1}A_{i\Gamma},\; i=2,3$ can be done in $O(N)$ operations instead of $O(N\log N)$. A similar situation arises when multiplying a vector by the matrix $A_{1\Gamma}^{\mathrm{T}}A_{11}^{-1}A_{1\Gamma}$. This is explained by the fact that when solving the problem $A_{11}\mathbf{\hat{y}}=A_{1\Gamma}\mathbf{\hat{f}}$ an essentially lesser number of iterations of the GMRES($k$) method is required as compared to $A_{11}\mathbf{\hat{y}}=\mathbf{\hat{f}}$. To solve the equation with the matrix $A_{11}$ 13 iterations of the  GMRES(k) method for the right-hand side $\mathbf{\hat{f}}$ were used, while for $A_{1\Gamma}\mathbf{\hat{f}}$, the number of iterations was $1$. Thus, each iteration of the CG method for solving problem (\ref{schur_decomp}) demands an essentially lesser number of arithmetical operations than one would use the Laplace operator as preconditioner for the whole computational domain. Moreover, in the latter case the number of iterations would be essentially larger due to a high contrast of the medium.

\begin{table}[!h]
\center \small
\begin{tabular}{llllllllll}

  \hline
    $N_z \times N_r$ & \multicolumn{3}{c}{$\begin{array}{l}3046\times 16384,\\d=17\end{array}$}& \multicolumn{3}{c}{$\begin{array}{l} 6078 \times 32768,\\ d=33\end{array}$}& \multicolumn{3}{c}{$\begin{array}{l}12134\times 65536,\\d=51\end{array}$}\\ \cmidrule(r){2-4} \cmidrule(r){5-7} \cmidrule(r){8-10}
    NP&$ {\mathrm{P}}$ & $\mathrm{T}$& $\mathrm{S}$&$ {\mathrm{P}}$ & $\mathrm{T}$& $\mathrm{S}$&$ {\mathrm{P}}$  & $\mathrm{T}$& $\mathrm{S}$\\
   32&42&9&-&179&38&-&-&-&-\\
   64&20&4.5&64&86&17.2&70&-&-&-\\
   128&20&2.1&137&60&9&135&379&53&-\\
   256&23 &1&288&46&4.5&270&197&27&251\\
   512& 34&0.65&443&56&2.1&579&383&13&521 \\
   1024&60 &0.56&514&83&1.52&800&281&8&848\\ \hline
  \end{tabular}
\caption{$\mathrm{\mathbf{P}}$ is the total time of the preliminaries, $\mathrm{\mathbf{T}}$ is the time of computing one harmonic from (\ref{series_lag1}), $\mathbf{S}$ is the speedup value, $\mathrm{\mathbf{d}}$ is the bandwidth of the precondition matrix,
$\mathrm{\mathbf{NP}}$ is the number of processors, $\mathrm{\mathbf{N_r,N_z}}$ is the number of mesh size towards $R$ and $Z$, respectively.}
\label{table44}
\end{table}

The Dichotomy Algorithm at all the stages of solving problem (\ref{acoustic_problem}) provides a high performance and scalability of the proposed parallel algorithm. This allows us to carry out engineering calculations (Fig.~\ref{pic:snapshot}) based on efficient algorithms with the use of thousands of processors. It should be noted that the most efficient and at the same time difficult for parallel realization  numerical methods are used.
\section{Conclusion}
In this paper the new parallel algorithm for solving SLAEs with the same block-tridiagonal matrix but different right-hand sides is proposed. To demonstrate the efficiency of the approach proposed, a problem of modelling the acoustic wave fields by the spectral-difference algorithm has been solved. A high performance of the Dichotomy Algorithm allows an effective use of the domain decomposition on a supercomputer. It should be noted that the domain decomposition was realized not for providing the parallel computation, but for decreasing the total number of arithmetical operations. In our case, the number of processors and subdomains are independent quantities, therefore the rate of convergence of the iterative method is independent of the number of processors. To solve the system of equations for the PML boundary conditions, the Dichotomy Algorithm was used.

To reduce the total computation time, the probing technique was used as preconditioning procedure. The probing technique in the context of parallel algorithms has not been widespread by now due to the necessity of solving SLAEs with band matrices. The efficient inversion of such matrices with the use of supercomputer systems is a non-trivial task. However, the development of the Dichotomy Algorithm has allowed one to overcome this difficulty. Now, this type of a preconditioner can be successfully implemented on supercomputers.

The numerical experiments carried out with $16$ up to $2048$ processors have proved the efficiency of the approach proposed. The dependence of the speedup value on the number of processors appears to be near-linear. Thus, a high performance and simplicity of service of the Dichotomy Algorithm allow one to include it into already existing sequential numerical procedures for their parallelization.

\newpage
\bibliography{base}

\begin{thebibliography}{10}

\bibitem{terekhov:Dichotomy}
A.~V. Terekhov.
\newblock Parallel dichotomy algorithm for solving tridiagonal system of linear
  equations with multiple right-hand sides.
\newblock {\em Parallel Comput.}, 36(8):423--438, 2010.

\bibitem{Samarski_Nikolaev}
A.A. Samarskij and E.S. Nikalayev.
\newblock {\em Numerical Methods for Grid Equations}.
\newblock Birkhauser Verlag, 1989.

\bibitem{Golub1989}
G.~H. Golub and C.~F. Van~Loan.
\newblock {\em Matrix computations (3rd ed.)}.
\newblock Johns Hopkins University Press, Baltimore, MD, USA, 1996.

\bibitem{Marsh1981}
F.~Marsh and D.~E. Potter.
\newblock Recurrence solution of a block tridiagonal matrix equation with
  neumann, dirichlet, mixed or periodic boundary conditions.
\newblock {\em Comp. Phys. Comm.}, 24:185--190, 1981.

\bibitem{block:Ruggieroa}
V.~Ruggieroa and E.~Galligania.
\newblock A parallel algorithm for solving block tridiagonal linear systems.
\newblock {\em Computers and Mathematics with application}, 24(4):15--21, 1992.

\bibitem{block:Akimova}
E.N. Akimova.
\newblock Parallel gauss algorithms for block tridiagonal linear systems.
\newblock {\em Matematicheskoe Modelirovanie}, 6(9):61--67, 1994.
\newblock (In Russian).

\bibitem{block:Hirshman}
S.P. Hirshman, K.S. Perumalla, V.E. Lynch, and R.~Sanchez.
\newblock Bcyclic: A parallel block tridiagonal matrix cyclic solver.
\newblock {\em J. Comp. Phys.}, 229:6392--6404, 2010.

\bibitem{block:Manyu}
X.~Manyu and L.~Quanyi.
\newblock A parallel iterative method for solving periodical block-tridiagonal
  linear equations.
\newblock {\em Applied Mathematics and Computation}, 184:599--607, 2007.

\bibitem{block:Mehrmann}
V.~Mehrmann.
\newblock Divide and conquer methods for block tridiagonal systems.
\newblock {\em Parallel Comput.}, 19(3):257--279, 1992.

\bibitem{Chung1995}
K.~L. Chung, Y.~H. Tsai, and W.~M. Yan.
\newblock A parallel solver for circulant block-tridiagonal systems.
\newblock {\em Computers \& Mathematics with Applications}, 29(1):109--113,
  1995.

\bibitem{Bai2007}
Y.~Bai and R.C. Ward.
\newblock A parallel symmetric block-tridiagonal divide-and-conquer algorithm.
\newblock {\em ACM Trans. Math. Softw.}, 33, August 2007.

\bibitem{Terekhov2}
A.~V. Terekhov.
\newblock Application of the parallel dichotomy algorithm for solving toeplitz
  tridiagonal systems of linear equations with one right-hand side(submitted).
\newblock {\em http://arxiv.org/abs/1002.2469}, 2010.

\bibitem{Terekhov2010}
A.~V. Terekhov.
\newblock High performance parallel algorithm for solving elliptic equations
  with non-separable variables.
\newblock {\em http://arxiv.org/abs/1002.3094v7 (Preprint)}, 2010.

\bibitem{fatab2011}
A.G. Fatyanov and A.V. Terekhov.
\newblock High-performance modeling acoustic and elastic waves using the
  parallel dichotomy algorithm.
\newblock {\em J. Comp. Phys.}, 230(5):1992--2003, 2011.

\bibitem{Tarek2008}
M.~Tarek.
\newblock {\em Domain Decomposition Methods for the Numerical Solution of
  Partial Differential Equations}, volume~61 of {\em Lecture Notes in
  Computational Science and Engineering}.
\newblock Springer Berlin Heidelberg, 2008.

\bibitem{Toselli:2004:DDM}
A.~Toselli and O.~Widlund.
\newblock {\em Domain Decomposition Methods - Algorithms and Theory}, volume~34
  of {\em Springer Series in Computational Mathematics}.
\newblock Springer, 2004.

\bibitem{Quarteroni:1999:DDM}
A.~Quarteroni and A.~Valli.
\newblock {\em Domain Decomposition Methods for Partial Differential
  Equations}.
\newblock Oxford Science Publications, 1999.

\bibitem{Konovalov}
N.~N. Yanenko, A.~N. Konovalov, A.~N. Bugrov, and G.~V. Shustov.
\newblock On organizing parallel computing and sweep parallelization.
\newblock {\em Chislennye Metody Mekhaniki Sploshnoi Sredy}, 9(7):139--146,
  1978.
\newblock (In Russian).

\bibitem{Wang2}
N.~Mattor, T.~J. Williams, and D.~W. Hewett.
\newblock Algorithm for solving tridiagonal matrix problems in parallel.
\newblock {\em Parallel Comput.}, 21(11):1769--1782, 1995.

\bibitem{MPI2}
M.~A. Heroux, P.~Raghavan, and H.~D. Simon.
\newblock {\em Parallel Processing for Scientific Computing (Software,
  Environments and Tools)}.
\newblock Society for Industrial and Applied Mathematics, Philadelphia, PA,
  USA, 2006.

\bibitem{top500}
http://www.top500.org/.

\bibitem{Smith:1996}
Barry~F. Smith, Petter~E. Bj{\o}rstad, and William Gropp.
\newblock {\em Domain Decomposition: Parallel Multilevel Methods for Elliptic
  Partial Differential Equations}.
\newblock Cambridge University Press, 1996.

\bibitem{Bitzarakis1997a}
S.~Bitzarakis, M.~Papadrakakis, and A.~Kotsopoulos.
\newblock Parallel solution techniques in computational structural mechanics.
\newblock {\em Comput. Methods Appl. Mech. Engrg.}, 148:75--104, 1997.

\bibitem{DD:GPU}
M.~Papadrakakis, G.~Stavroulakis, and A.~Karatarakis.
\newblock A new era in scientific computing: Domain decomposition methods in
  hybrid cpu–gpu architectures.
\newblock {\em Comput. Methods Appl. Mech. Engrg.}, 200:1490--1508, 2011.

\bibitem{Sheriff1995}
R.~E. Sheriff and L.~P. Geldart.
\newblock {\em Exploration Seismology}.
\newblock Cambridge University Press, 2nd edition, 1995.

\bibitem{Berenger1994}
J.P. Berenger.
\newblock A perfectly matched layer for the absorption of electromagnetic
  waves.
\newblock {\em J. Comp. Phys.}, 114:185--200, 1994.

\bibitem{Chew1994}
W.C. Chew and W.H. Weedon.
\newblock A 3d perfectly matched medium from modified maxwell's equations with
  stretched coordinates.
\newblock {\em Micro. Opt. Tech. Lett.}, 7:599--604, 1994.

\bibitem{Mikhailenko1999}
B.~G. Mikhailenko.
\newblock Spectral laguerre method for the approximate solution of time
  dependent problems.
\newblock {\em Applied Mathematics Letters}, 12:105--110, 1999.

\bibitem{abramowitz+stegun}
M.~Abramowitz and I.~A. Stegun.
\newblock {\em Handbook of Mathematical Functions with Formulas, Graphs, and
  Mathematical Tables}.
\newblock Dover, New York, ninth dover printing, tenth gpo printing edition,
  1964.

\bibitem{Samarskii2001}
A.A. Samarskii.
\newblock {\em The Theory of Difference Schemes}.
\newblock Marcel Dekker, 2001.

\bibitem{Strikwerda2004}
J.~C. Strikwerda.
\newblock {\em Finite Difference Schemes and Partial Differntial Equations}.
\newblock SIAM, 2 edition, 2004.

\bibitem{FEM2}
O.~C. Zienkiewicz and R.~L. Taylor.
\newblock {\em The finite element method}.
\newblock Butterworth-Heinemann, 2000.

\bibitem{Saad}
Y.~Saad.
\newblock {\em Iterative Methods for Sparse Linear Systems}.
\newblock SIAM, 2003.

\bibitem{Glowinski2007}
R.~Glowinski and Yu. Kuznetsov.
\newblock Distributed lagrange multipliers based on fictitious domain method
  for second order elliptic problems.
\newblock {\em Comput. Methods Appl. Mech. Engrg.}, 196:1498--1506, 2007.

\bibitem{Saulev1963}
V.~Saulev.
\newblock On solution of some boundary value problems on high performance
  computers by fictitious domain method.
\newblock {\em Siberian Math. J.}, 4:912--925, 1963.
\newblock (in Russian).

\bibitem{Ramiere}
I.~Ramiere, P.~Angot, and M.~Belliard.
\newblock A fictitious domain approach with spread interface for elliptic
  problems with general boundary conditions.
\newblock {\em Comput. Methods Appl. Mech. Engrg.}, 196:766--781, 2007.

\bibitem{Comformal_Mapping}
N.~Papamichael and N.~Stylianopoulos.
\newblock {\em Numerical Conformal Mapping: Domain Decomposition and the
  Mapping of Quadrilaterals}.
\newblock World Scientific Publishing Company, 2010.

\bibitem{Astrakhantsev1978}
G.~Astrakhantsev.
\newblock Method of fictitious domains for a second-order elliptic equation
  with natural boundary conditions.
\newblock {\em USSR Comput. Math. Math. Phys.}, 18:117--121, 1978.

\bibitem{Marchuk1982}
G.I. Marchuk.
\newblock {\em Methods of numerical mathematics}.
\newblock Springer-Verlag, 1982.

\bibitem{FFT}
U.~Schumann.
\newblock Fast fourier transforms for direct solution of poisson's equation
  with staggered boundary conditions.
\newblock {\em J. Comput. Phys.}, 75:123--137, 1988.

\bibitem{chan}
F.~C. Chan and E.~Keyes.
\newblock Interface precontitionings for domain-decomposed convection-diffusion
  operators.
\newblock In T.F. Chan, R~Glowinski, Periaux J., and Widlund O., editors, {\em
  Domain Decomposition Methods for Partial Differential Equations}. SIA, 1989.

\bibitem{Chan1992}
T.~F. Chan and T.~P. Mathew.
\newblock The interface probing technique in domain decomposition.
\newblock {\em SIAM J. Matrix Anal. Appl.}, 13:212--238, 1992.

\end{thebibliography}

\end{document}